\documentclass[12pt,leqno]{article}
\usepackage{amssymb,amsmath,amsthm}
\usepackage[all]{xy}
\usepackage[alphabetic,lite]{amsrefs}
\numberwithin{equation}{section}

\usepackage[top=2cm,bottom=2cm,left=3cm,right=3cm,marginparsep=0.3cm,marginparwidth=3cm,includefoot]{geometry}
\usepackage{graphicx}

\newcommand{\spa}{\vspace{0.5ex}\noindent}
\newcommand{\scbul}{\,\raise.4ex\hbox{$\scriptscriptstyle\bullet$}\,}
\def\shd{\mathcal{D}}
\def\shm{\mathcal{M}}
\def\shs{\mathcal{S}}
\newcommand{\Sol}{\shs ol}
\def\sho{\mathcal{O}}
\newcommand{\Ot}[1][X]{\sho^{\mspace{2.5mu}{\mathrm tp}}_{#1}}
\newcommand{\C}{{\mathbb{C}}}
\newcommand{\R}{{\mathbb{R}}}
\newcommand{\rhom}[1][]{{R\mathcal{H}om}_{\raise1.5ex\hbox to.1em{}#1}}
\newcommand{\hhom}{{\mathcal{H}om}}
\newcommand{\Thom}[1][]{{{\rm T}\mathcal{H}om}_{\raise1.5ex\hbox to.1em{}#1}}

\date{ICM 2018 Proceedings}
\begin{document}
\title{Fifty years of Mathematics with Masaki Kashiwara}
\author{Pierre Schapira}

\maketitle

Professor Masaki Kashiwara is certainly one of the foremost mathematicians of  our time. His influence is spreading over many fields of mathematics and the mathematical community slowly begins to appreciate the importance of the ideas and methods he has
introduced.

\subsubsection*{Mikio Sato}
Masaki Kashiwara was  a student of Mikio Sato and I will  begin with a few words about Sato (see~\cites{An07, Sc07} for a more detailed exposition). 
\begin{figure}\centerline{
\begin{tabular}{cc}
\includegraphics[scale=.1437]{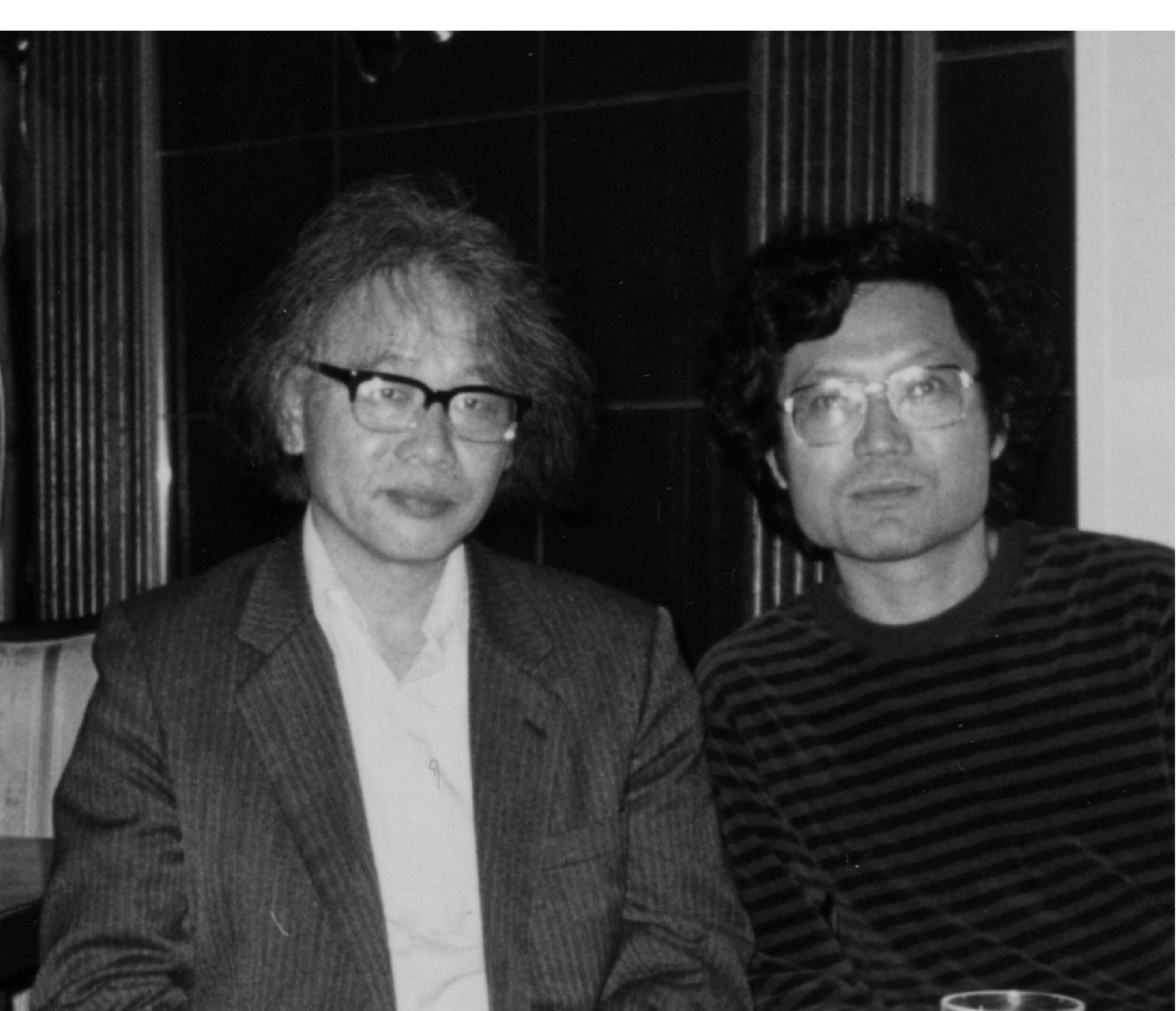} 
\end{tabular}}
\caption{Mikio Sato and Masaki Kashiwara}
\end{figure}
The story begins long ago, in the late fifties, 
when Sato created a new branch of mathematics
 now called  ``Algebraic Analysis'' by publishing 
 two papers on hyperfunction theory \cite{Sa60} and then 
developed his vision of analysis and linear partial differential
equations (LPDE) in a series of
lectures at the university of Tokyo in the 60s. Sato's idea is to define hyperfunctions on a real analytic manifold $M$ as 
 cohomology classes supported by $M$ 
of the sheaf $\sho_X$ of holomorphic functions on a complexification $X$ of $M$. 
One can then represent hyperfunctions as ``boundary values'' of holomorphic functions defined in tuboids in $X$ with wedge on $M$. To understand where the boundary values come from leads naturally Sato (see~\cite{Sa70}) to define his microlocalization functor and, as a byproduct,  the analytic wave front set of hyperfunctions. This is the starting point of microlocal analysis. Indeed, 
Lars H{\"o}rmander immediately understood the importance of Sato's ideas and adapted them to the $C^\infty$-setting by replacing 
boundary values of holomorphic functions with the Fourier transform (see~\cite{Ho71}).

In these old times, trying to understand real phenomena by complexifying a real manifold and looking at what happens in the complex domain was a totally  new idea.  And using cohomology of sheaves in analysis was definitly a revolutionary vision.

\subsubsection*{Master's thesis and the SKK paper}
Then came Masaki. In his master's thesis, dated 1970 and published in English in~\cite{Ka70},  he introduces and develops the theory of $\shd$-modules. Of course, a $\shd$-module is a module (right or left) over the non commutative sheaf of rings $\shd_X$ of holomorphic finite order differential operators on a given complex manifold $X$. And, as it is well known, a finitely presented module 
$M$ over a ring $R$ is the intrinsic way to formulate what is a finite system of $R$-linear equations with finitely many  unknowns. Hence, a coherent  $\shd_X$-module $\shm$ on $X$ is nothing but a system of linear partial  differential equations with holomorphic coefficients. Locally on $X$, it can be represented, non uniquely, by a matrix of differential operators. 

In this thesis, Masaki  defines
the operations of inverse or direct images for $\shd$-modules. Roughly speaking, these operations describe the system of equations satisfied by the restriction or the integral of the solutions of a system of equations.  
In particular he extends  the classical Cauchy-Kowalevski theorem to general
systems of LPDE. Consider the contravariant  functor $\Sol$, which to a $\shd_X$-module $\shm$ on the manifold $X$, associates the complex  (in the derived category of sheaves of $\C$-vector spaces) $\rhom[\shd_X](\shm,\sho_X)$ of its holomorphic solutions. The  Cauchy-Kowalevski-Kashiwara theorem  essentially asserts that the functor $\Sol$ commutes with the functor of taking the inverse image, under a non characteristic hypothesis. 

Hence Masaki Kashiwara may
be considered as the founder of analytic $\shd$-module theory, in parallel with Joseph Bernstein (see~\cite{Be71}) for the algebraic case, a theory which is now a fundamental tool in
many branches of mathematics, from number theory to mathematical physics. 

\medskip
The seventies are, for the analysts, the era of microlocal analysis. As mentioned above, the starting point was the introduction by Sato of  the microlocalization functor and the analytic wave front set. These ideas were then systematically developed in the famous paper~\cite{SKK73} by Mikio Sato, Takahiro Kawai and Masaki Kashiwara. Two fundamental results are proved here. 

\spa
First, the involutivity of characteristics of microdifferential systems. This was an open and fundamental question which had, at that time, only a partial answer due to  Quillen, Guillemin and Sternberg~\cite{GQS70}
(later, a purely algebraic proof was given by Gabber~\cite{Ga81}). 

\spa
The second result is a classification at generic points of any system of microdifferential equations. Roughly speaking, it is proved that, generically and after a so-called quantized contact transform, any such system is equivalent to a combination of 
a partial De Rham system, a partial Dolbeault system and a Hans Lewy's type system. 

This paper has had  an enormous influence on the analysis of partial differential equations (see in particular~\cites{Ho83, Sj82}).

\subsubsection*{The Riemann-Hilbert correspondence (regular case)}
Since the characteristic variety of a coherent $\shd$-module is involutive (one better say nowadays ``co-isotropic''), it is natural to look at the extreme case, when this variety is Lagrangian. One calls such systems ``holonomic''. They are the higher dimensional version of classical ordinary differential equations (ODE). Among ODE, there is a class of particular interest, called the class of Fuchsian equations, or also, the equations with regular singularities. Roughly speaking, the classical Riemann-Hilbert correspondence (R-H correspondence, for short) is based on the following question: given a finite set of points on the Riemann sphere and at each point, an invertible matrix of complex numbers (all of the same size), does there exist a unique Fuchsian ODE whose singularities are the given  points and such that the monodromy of its holomorphic solutions are the given matrices. 

From 1975 to 1980, Masaki Kashiwara gives a precise formulation of the conjecture establishing this correspondence in any dimension and eventually proves it. 

In 1975 (see~\cite{Ka75}) he proves that the contravariant  functor $\Sol$,  when restricted  to the derived category of $\shd$-modules with  holonomic cohomology, takes its values in the derived category of sheaves with $\C$-constructible cohomology. 
In the same paper, he also proves that if one starts with a $\shd$-module  ``concentrated in degree $0$'', then 
 the complex one obtains satisfies, what will be called five years later by Beilinson, Bernstein, Deligne and Gabber~\cite{BBD82},  the perversity conditions . 

Moreover, already in 1973, he gives in~\cite{Ka73} a formula to
calculate the local index of the
complex $\Sol(\shm)$ in terms of the characteristic cycle of the holonomic $\shd$-module $\shm$ and his formula contains the notion of ``local Euler obstruction'' introduced independently by MacPherson~\cite{McP74}. 

Classical examples in dimension $1$ show that the functor $\Sol$ cannot be fully faithful and the problem of defining the category of  ``regular holonomic $\shd$-modules'', the higher dimensional version of the Fuchsian ODE, remained open. 
For that purpose, Masaki introduces with Toshio Oshima~\cite{KO77} the notion of regular singularities along a smooth involutive manifold and then formulates  precisely in 1978 the notion of regular holonomic $\shd$-module and  what should be the  Riemann-Hilbert correspondence (see~\cite{Ra78}), namely an equivalence of categories between the derived category of $\shd$-modules with regular holonomic cohomology and the derived category of sheaves with $\C$-constructible cohomology. 
He  solves this conjecture in 1980 (see~\cites{Ka80, Ka84}) by constructing a quasi-inverse to the functor $\Sol$,  
 the functor $\Thom$ of tempered cohomology. For a constructible sheaf $F$, the object   $\Thom(F,\sho_X)$ is represented by applying the functor $\hhom(F,\scbul)$ to the Dolbeault resolution of $\sho_X$ by differential forms with distributions as coefficients. 

Of course,  Kashiwara's paper came after  Pierre Deligne's famous book~\cite{De70} in which he solves the R-H problem for regular connections. This book  has had a deep influence on  the  microlocal approach of the R-H correspondence,  elaborated by Masaki jointly with  T.~Kawai~\cite{KK81}. Finally note that a different proof of this theorem was obtained later by Zogman Mebkhout in~\cite{Me84}.

\subsubsection*{Other results on $\shd$-modules and related topics}
\begin{itemize}
\item[{\rm (i)}]
Besides his proof of the R-H correspondence, Masaki obtains fundamental results in $\shd$-module theory. He proves in~\cite{Ka76}  the rationality of the zeroes of the $b$-function of Bernstein-Sato by using Hironaka's theorem and adapting Grauert's direct image theorem to $\shd$-modules. 
\item[{\rm (ii)}]
Motivated   by the theory of holonomic $\shd$-modules,
Masaki proves in \cite{Ka82}  the codimension-one property of
quasi-unipotent sheaves. 
\item[{\rm (iii)}]
Masaki gives a fundamental contribution to the theory of
 ``variation of (mixed) Hodge structures''  (see for example~\cites{Ka85b, Ka86}).
\item[{\rm (iv)}]
In~\cite{BK86}, Masaki and Daniel Barlet endow  regular holonomic $\shd$-modules with a ``canonical'' good filtration. 
\item[{\rm (v)}]
A classical theorem of complex geometry (Frisch-Guenot, Siu, Trautmann) asserts that,  on a complex manifold $X$, any reflexive coherent sheaf  defined on the complementary of a complex subvariety of codimension  at least $3$ extends as a coherent sheaf through this subvariety.  
The codimension $3$ conjecture is an analogue statement for holonomic microdifferential modules when replacing $X$ with a Lagrangian subvariety of the cotangent bundle. This extremely difficult conjecture was recently proved by Masaki together with Kari Vilonen in~\cite{KVi14}.
\item[{\rm (vi)}]
Kashiwara's book on $\shd$-modules~\cite{Ka03} contains a lot of original and deep results. In this book he defines in particular the microlocal $b$-functions and gives a tool, the  ``holonomy diagrams'', to calculate them. 
\item[{\rm (vii)}]
The book on category theory~\cite{KS06}, written with P.~Schapira, sheds  new light on a  very classical subject and contains a great deal of original results. 
\end{itemize}

\subsubsection*{Mathematical physics} 
\begin{itemize}
\item[{\rm (i)}]
In collaboration with Takahiro~Kawai and Henri~Stapp, Masaki applied the theory of holonomic $\shd$-modules to the study 
of Feynman integrals. See in particular~\cites{KK77, KK77b, KK78, KKS77}.
\item[{\rm (ii)}]
In~\cite{Sa82}, Mikio Sato and Yasuko Sato established that soliton equations are dynamical systems on the infinite Grassmann manifold. Based on this work, Kashiwara, with Etsuro~Date, Michio~Jimbo and Tetsuji~Miwa (see~\cites{DJKM81, DJKMb81, DJKM82}), have found links between hierarchies of soliton equations and representations of infinite dimensional Lie algebras, {\em e.g.,} between the KP hierarchy and $\mathfrak{gl}_\infty$, the KdV hierarchy and the affine Lie algebra of type $A_1^{(1)}$, and so on. In terms of Hirota's dependent variable \cite{Hi71}, the set of soliton solutions of a hierarchy is identified with the group orbit of $1$ in the space of the vertex operator representation of the corresponding infinite-dimensional Lie algebra\footnote{I warmly thank Tetsuji Miwa for his help concerning this topic.}. 
\end{itemize}

\subsubsection*{Representation theory}
 Masaki Kashiwara also had an enormous influence in representation theory, harmonic analysis and quantum groups. His work has transformed the field, in its algebraic, categorical, combinatorial, geometrical and analytical aspects.
\begin{itemize}
 \item[{\rm (i)}]
In~\cite{KKMOO78}, Kashiwara solves a conjecture of Helgason on non-commutative harmonic analysis.
\item[{\rm (ii)}]
At the same period, he proves a fundamental result on the Campbell-Hausdorff formula~\cite{KV78}  in collaboration with Mich{\`e}le Vergne. There is currently a lot of activity stemming from this paper.
\item[{\rm (iii)}]
In collaboration with Jean-Luc Brylinski, he solves in~\cite{BK81} a major open problem in representation theory, the Kazhdan-Lusztig conjecture on infinite-dimensional representations of simple Lie algebras, a conjecture proved independently by Beilinson-Bernstein in~\cite{BB81}. This is one of the most influential paper in geometric representation theory.
\item[{\rm (iv)}]
These results are generalized
to Kac-Moody algebras with Toshiyuki Tanisaki \cite{KT90}: this 
was one of the key steps in the proof of Lusztig's conjecture on simple modules
for algebraic groups in positive characteristic.
\item[{\rm (v)}]
Kashiwara has also obtained major results on representations of real Lie groups.

He reinterprets  the Harish-Chandra theory in terms of $\shd$-module theory and
obtains by this method important theorems on semi-simple Lie algebras
with Ryoshi Hotta~\cite{HK84},  on real reductive groups with  Wilfried Schmid~\cite{KSw94}.

The final theory constructed by Kashiwara  \cite{Ka08} shows how to construct geometrically
the Lie group representations coming from Harish-Chandra modules: the first step
is localization, turning Harish-Chandra modules into
$\shd$-modules on the flag variety. Kashiwara's Riemann-Hilbert correspondence turns those into
constructible sheaves. Via a sheaf  theoretic version of the Matsuki correspondence, these become
equivariant sheaves for the real Lie group, which lead to the correct 
representations of the real Lie group.
\item[{\rm (vi)}]
In \cite{KR08}, using deformation quantization modules (see below), Kashiwara constructs with Rapha{\"e}l Rouquier  a
microlocalization of rational Cherednik algebras. This is the first extension of classical localization methods to symplectic manifolds that are not cotangent bundles and  opens a new direction in geometric representation theory.  
\end{itemize}

\subsubsection*{Quantum groups and crystal bases\footnote{I warmly thank Rapha{\"e}l Rouquier for his help on this section.}}

\begin{itemize}
\item[{\rm (i)}]
Finite-dimensional
representations of compact Lie groups are some of the most fundamental objects
in representation theory. The search for good bases in these representations,
in relation with invariant theory and geometry, was a source of attention
since the late 19th century. A change of paradigm occurred with 
Kashiwara's work in 1990. This is based on quantum groups, which are
deformations of the enveloping algebras of Kac-Moody Lie algebras.
Kashiwara discovered that, when the quantum group parameter goes to $0$ (temperature zero limit in the  solvable lattice models setting of statistical mechanics), the theory acquires a combinatorial structure,
replacing the linear structure. That leads to a basis at parameter $0$ (crystal basis), whose existence was proven by an extraordinary combinatorial tour-de-force~\cite{Ka90}.
Note that George Lustzig also considered the bases at $q=0$ given by the PBW
bases in the ADE case, and constructed canonical bases (see~\cite{Lu90}).

Crystal  basis can be lifted uniquely to a basis (global basis) satisfying certain 
 symmetry properties. 
Crystal bases are now a central chapter of representation theory and algebraic combinatorics. 

\item[{\rm (ii)}]
Kashiwara used those crystal bases to solve a very basic problem of representation theory, the decomposition of tensor products of irreducible representations of simple Lie algebras.

\item[{\rm (iii)}]
Kashiwara has given with Yoshihisa Saito a geometric
construction of the crystal basis in terms of Lagrangian
subvarieties of  Lusztig  quiver varieties \cite{KSa97}.

\item[{\rm (iv)}]
Kashiwara's recent work on higher representation theory
has been fundamental.  In~\cite{KK12},  he solves with Seok-Jin  Kang
a basic open problem: he proves
that cyclotomic quiver Hecke algebras give rise to simple 2-representations of Kac-Moody
algebras.

\item[{\rm (v)}]
Kashiwara has obtained a number of key results on finite-dimensional
representations of affine quantum groups, in particular on the irreducibility
of tensor products. This has led to new directions in higher representation theory.

One such instance is the groundbreaking discovery by Kashiwara, together with Seok-Jin  Kang and Myungho  Kim (see~\cite{KKK18}) of a new type of Schur-Weyl duality relating quantum affine algebras of arbitrary
types and certain quiver Hecke algebras. Another one is the general construction of 
monoidal categorification of cluster algebras via quiver Hecke algebras with these two authors and Se-jin Oh (\cite{KKKO18}).
\end{itemize}

\subsubsection*{Microlocal sheaf theory}
From 1982 to 1990, with Pierre Schapira,  he introduces and develops   the microlocal theory of sheaves (see~\cites{KS82, KS85, KS90}). This theory emerged from a joint paper (see~\cite{KS79}) in which they solve the Cauchy problem for microfunction solutions of hyperbolic $\shd$-modules on a real analytic manifold. Indeed, the basic idea is that of microsupport of sheaves which gives a precise meaning to the concept of propagation. On a real manifold $M$, for  a (derived) sheaf $F$, its microsupport, or  singular support, is a closed conic subset of the cotangent bundle $T^*M$ which describes the codirections of non extension of sections of $F$. The microsupport of sheaves is, in some sense, a real analogue of the characteristic variety of coherent $\shd$-modules on complex manifolds and the functorial properties of the microsupport are very similar to those of the characteristic variety of $\shd$-modules. The precise link between both notions is a result which asserts that the microsupport of the complex $\Sol(\shm)$ of holomorphic solutions of a coherent $\shd$-module $\shm$ is nothing but the characteristic variety of $\shm$. 
Moreover, and this is one of the main results of the theory,  
the microsupport  is co-isotropic. As a by-product, one obtains  a completely different proof of the involutivity of characteristics of $\shd$-modules. 
\begin{figure}\centerline{
\begin{tabular}{cc}
\includegraphics[scale=.1437]{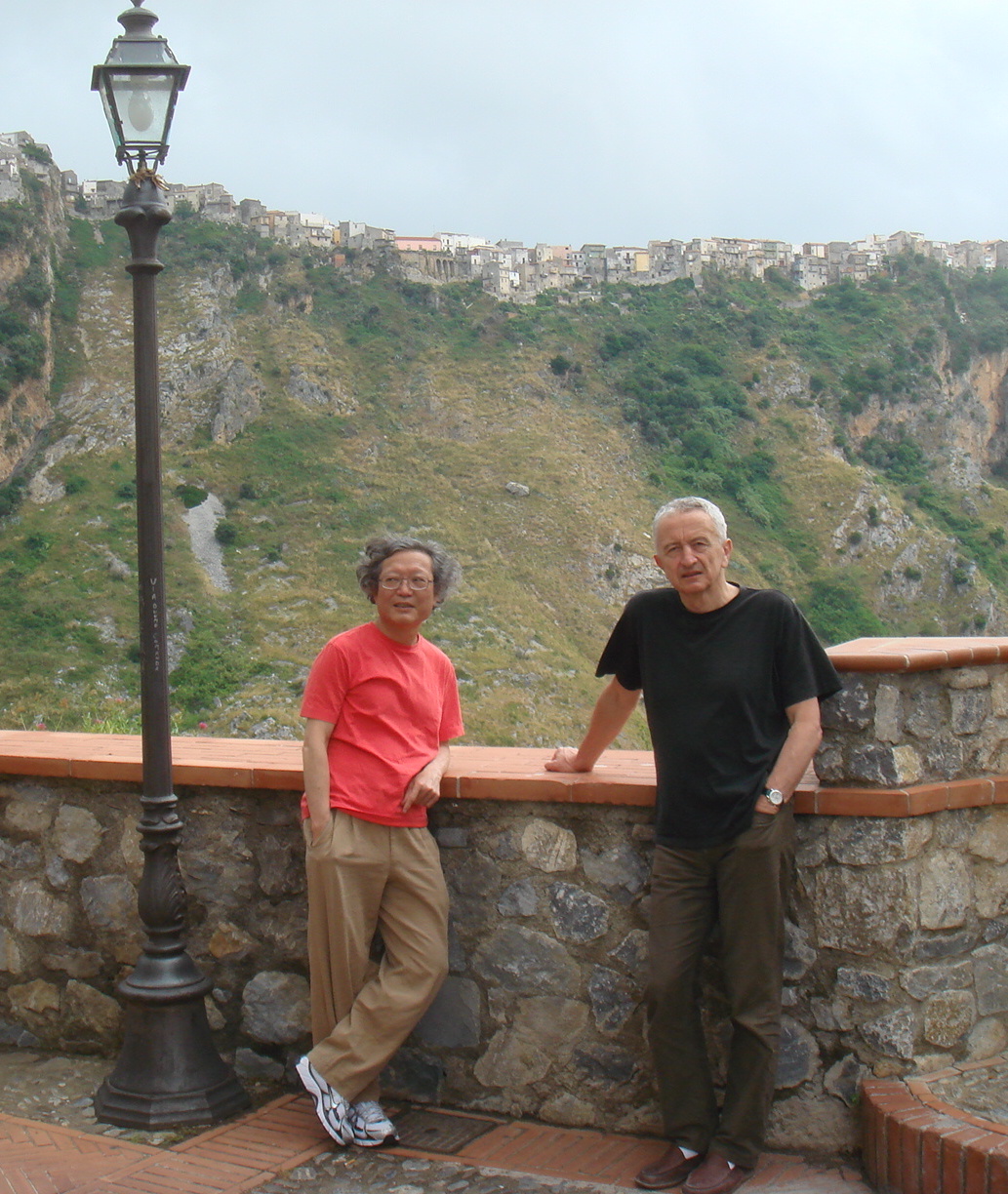} 
\end{tabular}}
\caption{Masaki Kashiwara and Pierre Schapira, Italy 2009}
\end{figure}

By using the microsupport, one can now   localize the derived category of sheaves on open subsets of $T^*M$ and the prestack (presheaf of categories) one obtains is a candidate to be a first step for an alternative construction of the Fukaya category, a program recently initiated by Dmitry Tamarkin~\cite{Ta15}. 

Microlocal sheaf theory naturally leads Kashiwara to extend his previous work on complex analytic Lagrangian cycles to the real setting. In~\cite{Ka85} he defines  the characteristic cycle of an $\R$-constructible
sheaf and gives a new index
formula. He also gives in this context a remarkable and unexpected
``local Lefschetz formula'' with
applications to representation theory (see~\cite{KS90}*{Ch.~IX~\S~6}).

Microlocal sheaf theory has found applications in many other fields of mathematics, such as representation theory (see above) and  symplectic topology with Tamarkin, Nadler, Zaslow, Guillermou  and many others (see in particular~\cites{Ta08, NZ09, Na09, Gu12}).  It has also applications  in knot theory thanks to a result of~\cite{GKS12} which implies that the category of simple sheaves along a smooth Lagrangian submanifold is a Hamiltonian isotopy invariant (see {\em e.g.,}~\cite{STZ17}). Recently, Alexander Beilinson~\cite{Be16} adapted the definition of the microsupport of sheaves to  
arithmetic geometry and the theory is currently  being developed, in particular by Takeshi Saito (see~\cite{Sai17}).

\subsubsection*{Deformation quantization}
In~\cite{Ka96}, Masaki  introduces  the notion of algebroid stacks in order to 
quantize complex contact manifolds several years before such
constructions become extremely popular. 

Then, with  P.~Schapira, they undertook in the book~\cite{KS12}  a systematic study of DQ-modules (DQ for 
``Deformation Quantization'')  on complex Poisson manifolds. 
This is a theory which contains both that of usual $\shd$-modules  and classical analytic geometry (the commutative case). A perversity theorem  (in the symplectic case) is obtained. This book also contains the precise statement of an old  important conjecture of Masaki on the Todd class in the Riemann-Roch theorem, a conjecture recently proved by Julien Grivaux~\cite{Gr12}  (see also Ajay Ramadoss~\cite{Ra08} for the algebraic case).

An illustration of the usefulness  of DQ-modules is the quantization of Hilbert schemes of points on the plane,  constructed
in~\cite{KR08}. 

\subsubsection*{Ind-sheaves and the irregular Riemann-Hilbert correspondence}
As already mentioned,  Masaki introduced the functor of
tempered  cohomology in the 80s, in his proof of the R-H correspondence. This functor is systematically studied with P.~Schapira in~\cite{KS96} where a dual functor, the functor of Whitney tensor product, is 
also  introduced. However, the construction of these two functors appears soon as a particular case  of a more general notion, that of ind-sheaves, that is, ind-objects of the category of sheaves with compact supports. This theory  is developed
 in~\cite{KS01}. 

Ind-sheaf theory  is a tool to treat functions or distributions  with growth conditions with the techniques of sheaf theory. In particular, it allows one to define the (derived) sheaf  $\Ot[X]$ of holomorphic functions with tempered growth (a sheaf for the so-called subanalytic topology).  Already, in the early 2000, it became clear that this ind-sheaf  was an essential tool for  the study of irregular holonomic modules and a toy model was studied  in~\cite{KS03}. However, although the functor of tempered holomorphic solutions is much more precise than the usual functor $\Sol$, it is still not enough precise to be fully faithful on the category of irregular holonomic $\shd$-modules. Then, by adapting to ind-sheaves a construction of Dmitry Tamarkin~\cite{Ta08},  Masaki and Andrea D'Agnolo introduced in~\cite{DK16} the ``enhanced ind-sheaf of tempered holomorphic functions'' and obtained a fully faithful functor. 
This deep theory, which uses in an essential manner the  fundamental results of Takuro Mochizuki~\cites{Mo09, Mo11} (see also Claude Sabbah~\cite{Sa00} for preliminary results  and Kiran Kedlaya~\cites{Ke10,Ke11} for the analytic case), has important applications, in particular in the study of the Laplace transform (see~\cite{KS16}). 
\begin{figure}\centerline{
\begin{tabular}{cc}
\includegraphics[scale=.1437]{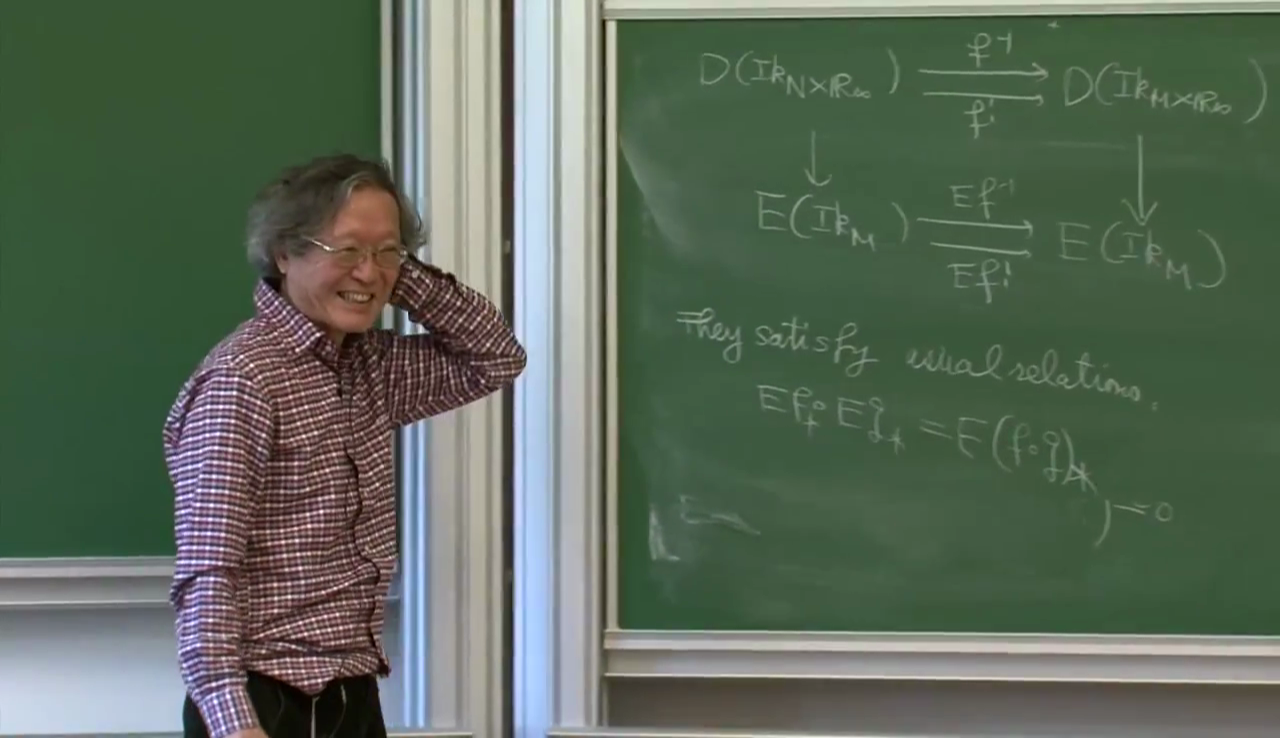} 
\end{tabular}}
\caption{Masaki Kashiwara, IHES 2015}
\end{figure}

\subsubsection*{Conclusion}
Kashiwara's contribution to mathematics is really astonishing and 
it should be mentioned that 
his influence is not only due to his published work, but also
to many informal talks. Important subjects such as 
second microlocalization, complex quantized contact transformations,
the famous``watermelon theorem'', etc. were initiated by him, although
not published.  Masaki is an invaluable source of
inspiration for many people.

\providecommand{\bysame}{\leavevmode\hbox to3em{\hrulefill}\thinspace}

\vspace*{1cm}
\noindent
\parbox[t]{21em}
{\scriptsize{
\noindent
Pierre Schapira\\
Sorbonne Universit{\'e}\\
Campus Pierre et Marie Curie, IMJ-PRG\\
4 place Jussieu, 75252 Paris Cedex 05 France\\
e-mail: pierre.schapira@imj-prg.fr\\
http://webusers.imj-prg.fr/\textasciitilde pierre.schapira/
}}

\end{document}